\documentclass[11pt]{article}
\usepackage{amsfonts, amsmath, amssymb, amsthm, algorithm,algpseudocode, enumitem}

\PassOptionsToPackage{obeyspaces}{url}
\usepackage[colorlinks=true,citecolor=blue,urlcolor=blue,linkcolor=blue,bookmarksopen=true]{hyperref}

\DeclareMathOperator{\ex}{ex}

\author{Zhiyang He\footnote{Department of Mathematical Sciences, Carnegie Mellon University, Pittsburgh, United States. Advised by Boris Bukh, supported by NSF grant DMS-1555149.}}
\title{New Upper Bound on Extremal Number of Even Cycles}

\newtheorem*{redLemma}{Reduction Lemma}
\newtheorem{theorem}{Theorem}
\newtheorem{claim}{Claim}
\newtheorem{definition}{Definition}
\newtheorem{remark}{Remark}
\newtheorem{lemma}[theorem]{Lemma}
\newtheorem{corollary}[theorem]{Corollary}
\floatname{algorithm}{Procedure}
                           
\newcommand*{\abs}[1]{\lvert #1\rvert}
\newcommand*{\eqdef}{\stackrel{\text{\tiny{def}}}{=}}
\newcommand*{\larrow}{\leftrightsquigarrow}

\begin{document}
\maketitle

\begin{abstract}
In this paper, we prove $\ex(n, C_{2k})\le (16\sqrt{5}\sqrt{k\log k} + o(1))\cdot n^{1+1/k}$. We improved on Bukh--Jiang's method used in their 2017 publication, thereby reducing the best known upper bound by a factor of $\sqrt{5\log k}$. 
\end{abstract}

\section*{Introduction}
In the field of extremal graph theory, Tur\'an's Problem, introduced by Tur\'an~\cite{Turan41} in 1941, asks the following question: Given graph $F$, what is the maximum number of edges that a graph on $n$ vertices can have while not containing $F$ as a subgraph? This number, denoted $\ex(n, F)$, is now referred to as the Tur\'an number or the extremal number of $F$. Similarly, for a family of graphs $\mathcal{F}$, $\ex(n, F)$ requires that no element of $\mathcal{F}$ is present.

The first result, known as Mantel's Theorem, was proven by Mantel~\cite{Mantel1907} in 1907. Since then, extensive amount of works have been established, among which is the celebrated Erd\H{o}s-Stone-Simonovits Theorem~\cite{ErdosStone}: $\ex(n, F) = (1-\frac{1}{\chi(F)-1}+o(1))\binom{n}{2}$. This result, proven in 1946, essentially solved Tur\'an's Problem for all graphs $F$ with $\chi(F) > 2$. However, the case for bipartite graphs is left open. To date, no results comparable to Erd\H{o}s-Stone-Simonovits Theorem have been derived. For the two most studied families of bipartite graphs, complete bipartite graphs $K_{s, t}$ and even cycles $C_{2k}$, the magnitude of $\ex(n, F)$ is not known for general $s, t, $ and $k$. In particular, it has been known for decades that $\ex(n, C_{2k}) = O_k(n^{1+1/k})$, while a matching lower bound has not been established. 

To discuss methods that lead to upper bounds on $\ex(n, C_{2k})$, we first show a simple derivation of $\ex(n, \{C_3, C_4, \cdots, C_{2k}\}) \le cn^{1+1/k}$ for some constant $c$. Consider a graph containing $\Theta(n^{1+1/k})$ edges with girth at least $2k+1$, and reduce it to a graph with minimum degree $\Theta(n^{1/k})$. Fix arbitrary vertex $v$, we start a Breadth-First Search (BFS) at $v$ and observe that for the first $k$ levels of the breadth-fist search tree, every level must expand by a factor of $\Theta(n^{1/k})$ compared to the previous level. In particular, no two vertices with depth less than $k$ can have common neighbors with greater depth. Since the $k$th level cannot have more than $n$ vertices, the bound follows. We present this derivation since the best upper bounds on $\ex(n, C_{2k})$ are, in essence, all established using this same approach. As we will see shortly, employing this method imposes fundamental limitations to the results derivable. 

The first important upper bound on $\ex(n, C_{2k})$ was proved by Bondy-Simonovits~\cite{BondySimonovits} in 1974, where they showed $\ex(n, C_{2k}) \le 20kn^{1+1/k}$. This result is subsequently improved through a line of researches, most recently by Pikhurko~\cite{Pikhurko} to $\ex(n, C_{2k}) \le (k-1)n^{1+1/k} + O_k(n)$ in 2010 and by Bukh--Jiang~\cite{BukhJiang} to $\ex(n, C_{2k}) \le 80\sqrt{k}\log kn^{1+1/k} + O_k(n)$ in 2017. Our main contribution in this paper is the following theorem. 
\begin{theorem}\label{thm:main}
Fix $k$, let $G$ be a $n$-vertex graph where $n\ge (20k)^{4k^3+2k^2}$. If 
\[
  \abs{E(G)} > 16\sqrt{5}\sqrt{k\log_e k}\cdot n^{1+1/k}+8000k^{4}n^{1+(2k-1)/(2k^2)},
\]
then $G$ contains a copy of $C_{2k}$. 
\end{theorem}

For the rest of this paper, we will abbreviate $\log_e$ as $\log$. Our approach is an improved version of Bukh--Jiang's approach, and therefore suffers the same limitation as all BFS arguments. More specifically, consider a bipartite graph $G$ with bipartition $V_1, V_2$ such that $\abs{V_1} = n, \abs{V_2} = n/(k-1)$. The BFS argument for the girth problem can be exploited to show that $e(G)\le c(k-1)^{-1/2}n^{1+1/k}$ if $G$ has girth at least $2k+1$ (For a more detailed argument, see \cite{NaorVerstraete2005}). Now if we duplicate each vertex in $V_2$ into $k-1$ copies, we obtain a graph on $2n$ vertices and $c\sqrt{k-1}n^{1+1/k}$ edges with no $C_{2k}$. Therefore, the best upper bound on $\ex(n, C_{2k})$ derivable from the BFS argument is $c\sqrt{k}n^{1+1/k}$ for fixed constant $c$. To break the $O(\sqrt{k}n^{1+1/k})$ threshold would require a different approach.

Our result improves the best known bound for $\ex(n, C_{2k})$ by a factor of $\sqrt{5\log k}$, taking us one step closer to the limitation of the method. Before discussing the proof, we would like to point out the following facts about this paper. This paper is modified from Bukh--Jiang's manuscript. While we had made global modifications to Bukh--Jiang's methods, most improvements made are local. In particular, it is highly similar to their paper mathematically, with a few statements and minor proofs largely unmodified. The author's intentions in writing this paper this way are to give a more intuitive delivery of Bukh--Jiang's methods, to present simplifications and improvements that lead to a better result, and to avoid confusing readers with different notations and proof structures that demonstrate the same ideas. Therefore, this paper adopted the same notation with some unmodified definitions from Bukh--Jiang's paper, and reshaped the delivery structure and language to uncover the underlying ideas and intuition.  

To begin our proof, in Section~\ref{sec:def} we describe the graph structures used in this paper. For a detailed discussion of how our methods relate and differ from Pikhurko's and Bukh--Jiang's work, please see Section~\ref{sec:preface}.

\section{Graph Reduction and Exploration}\label{sec:def}
To employ the breadth-first search approach, we first process our graph to gain control over the degrees of vertices. Classically, the graph is reduced to have minimum degree $O_k(n^{1/k})$ at the expense of half of the edges. However, in Bukh--Jiang's approach and in our approach, control over maximum degree is also required. Bukh--Jiang modified the BFS process to avoid vertices of high degrees, while we make use of the following reduction lemma.
\begin{redLemma}
Fix $\alpha\in (0, 1)$, let $\gamma=(20/\alpha)^{-2/\alpha}$. Let $d_{\text{min}}(G), d_{\text{max}}(G)$ denote the minimum and maximum degree of a graph $G$, respectively. If a graph $G$ on $n$ vertices has at least $cn^{1+\alpha}$ edges, then it contains a subgraph $G'$ such that $\abs{V(G')} \ge c\gamma n^{\alpha/2}$, $\abs{E(G)} \ge (c/4)v(G')^{1+\alpha}$, $d_{\text{min}}(G') \ge (c/2)v(G')^{\alpha}$, and $d_{\text{max}}(G')/d_{\text{min}}(G') \le 1/\gamma$. 
\end{redLemma}
An initial version of this lemma was first proved by Erd\"{o}s-Simonovits~\cite{ErdosSimonovits}, and various forms of this lemma occur in other works. Bukh--Jiang proved a slightly different version of this lemma in their addendum. By slightly modifying their proof, we obtain the above lemma. This proof is included in the Appendix for completeness. 

With this structure in mind, our real result in this paper is the following theorem.
\begin{theorem}\label{thm:real}
Fix $k\ge 4$, let $\Delta = \sqrt{k}(20k)^{2k}$, and let $d\ge \max(2\sqrt{5}\sqrt{k\log k} n^{1/k},\: (20k)^{4k^2+2k})$. If $G$ is a graph on $n$ vertices such that $d_{\text{min}}(G)\ge 2d + 5k^2$ and $d_{\text{max}}(G)\le \Delta d$, then $G$ contains a copy of $C_{2k}$. 
\end{theorem}

Theorem~\ref{thm:main} then follows from Reduction Lemma and Theorem~\ref{thm:real}.
\paragraph{Proof of Theorem~\ref{thm:main}}
Assume a graph $H$ on $m$ vertices has more than $16\sqrt{10}\sqrt{k\log k}\cdot n^{1+1/k}+8000k^{4}n^{1+(2k-1)/(2k^2)}$ edges, then we can find a bipartite subgraph $H'$ with at least half of its edges. Using the Reduction Lemma, we find a subgraph $G$ on $n\ge 4\sqrt{10}\sqrt{k\log k}\gamma m^{1/(2k)}$ vertices and at least $2\sqrt{10}\sqrt{k\log k}n^{1+1/k}$ edges, where $\gamma = (20k)^{-2k}$. Now we compute the minimum degree in $G$.

Let $c = e(H)/m^{1+1/k}$, we have that $c \ge 8\sqrt{10}\sqrt{k\log k} + 4000k^4/m^{1/(2k^2)}$, which implies that
\begin{align*}
    d_{\text{min}}(G) \ge \frac{c}{2}n^{1/k} \ge 4\sqrt{10}\sqrt{k\log k}n^{1/k} + \frac{2000k^4(\gamma m^{1/(2k)})^{1/k}}{m^{1/(2k^2)}} \ge 4\sqrt{10}\sqrt{k\log k}n^{1/k} + 5k^2.
\end{align*}
Now from Pikhurko's Result \cite{Pikhurko}, we know that if  $d_{\text{min}}(G) \ge kn^{1/k}$, then $G$ contains a $C_{2k}$. Therefore, Reduction Lemma implies $d_{\text{max}}(G) \le (20k)^{2k}kn^{1/k}$. Let $d = 2\sqrt{10}\sqrt{k\log k}n^{1/k}$, $\Delta = \sqrt{k}(20k)^{2k}$. Theorem~\ref{thm:real} completes the proof.
\qed

To prove Theorem~\ref{thm:real}, we elaborate further on the graph structures. Let $G$ be a graph as in the statement of Theorem~\ref{thm:real}. Fix arbitrary vertex $v$ of $G$ and start a breadth-first search process at $v$. Let $V_i$ be the set of vertices at minimum distance $i$ from $v$ for $i\in [k]$. We recall the following definition of a tri-layered graph, which is the basis of our discussions in Section~\ref{sec:theta}, from Bukh--Jiang. 
\begin{definition}[Bukh--Jiang~\cite{BukhJiang}]
A graph $G$ is called trilayered if its vertex set can be partitioned into $V_1, V_2, V_3$ such that all edges in $G$ are between $V_1, V_2$ or between $V_2, V_3$. For arbitrary $G$, we use $G[V_1, V_2, V_3]$ to denote the induced trilayered graph of $G$ on $V_1, V_2$ and $V_3$. For $A, B, C, D\in \mathbb{R}$, we say that a trilayered graph has minimum degree $[A : B, C : D]$ if the minimum degree from $V_1$ to $V_2$, $V_2$ to $V_1$, $V_2$ to $V_3$ and $V_3$ to $V_2$ are at least $A, B, C, D$, respectively. 
\end{definition}
The last ingredient we need is the following definition of a $\Theta$-graph, which is at the core of all our future discussions.
\begin{definition}
A $\Theta$-graph is a cycle of length at least $2k$ with a chord. That is, an edge outside of the cycle connecting two vertices of the cycle. 
\end{definition}
The rest of the paper is organized as follows: In Section~\ref{sec:preface}, we recall several important results from Pikhurko and Bukh--Jiang, which prove the non-existence of $\Theta$-graphs in the trilayered subgraphs formed by our breadth-first search exploration. In Section~\ref{sec:theta}, which contains our main improvements in this paper, we argue that if certain conditions hold, then a trilayered graph satisfies certain minimum degree condition must be present. We then embed a $\Theta$-graph in such subgraphs, contradicting our result from Section~\ref{sec:preface}. In Section~\ref{sec:final}, we show that either the aforementioned conditions hold, or the levels from exploration expand exponentially. Final computations then prove Theorem~\ref{thm:real}. 

\section{Results on \texorpdfstring{$\Theta$}{Theta}-Graphs}\label{sec:preface}
To argue for non-existence of $\Theta$-graphs in our exploration, we recall results of Pikhurko.
\begin{lemma}[Lemma~2.2 in~\cite{Pikhurko}]\label{lem:min_embed}
Let $k\ge 3$. Any bipartite graph $H$ of minimum degree at least $k$ contains a $\Theta$-graph.
\end{lemma}
\begin{corollary}\label{cor:avg_embed}
Let $k\ge 3$. Any bipartite graph $H$ of average degree at least $2k$ contains a $\Theta$-graph.
\end{corollary}
\begin{lemma}[Claim~3.1 in~\cite{Pikhurko}]\label{lem:pikhurko_notheta}
Suppose $G$ contains no $C_{2k}$. For $1\le i\le k-1$, neither of $G[V_i]$ and $G[V_i,V_{i+1}]$ contains a bipartite $\Theta$-graph.
\end{lemma}

Using these results, Pikhurko showed that every level must expand by a factor of roughly $d/k$ compared to the previous level. The bound $ex(n, C_{2k})\le O(kn^{1+1/k})$ then followed. Bukh and Jiang improved on his method by analyzing three consecutive levels, proving a better expansion ratio among them. They employed the following technical definition, which generalized $\Theta$-graphs to three levels, and proved the next lemma in conjunction.

\begin{definition}
Let $G$ be a trilayered graph with layers $V_1,V_2,V_3$. A $\Theta$-graph $T$ in $G$ is well-placed if every vertex of $T$ in $V_2$ is adjacent to some vertex of $V_1$ not in $T$. 
\end{definition}
\begin{lemma}[Lemma~10 in \cite{BukhJiang}]\label{lem:trilayer_notheta}
Suppose $G$ contains no $C_{2k}$. For $1\le i\le k-1$, the graph $G[V_{i-1},V_i,V_{i+1}]$ contains no well-placed $\Theta$-graphs.
\end{lemma}

Note that Lemma~\ref{lem:trilayer_notheta} is analogous to Lemma~\ref{lem:pikhurko_notheta}. To prove statements equivalent to Lemma~\ref{lem:min_embed} in trilayered graphs,
Bukh--Jiang analyzed trilayered subgraphs with specific minimum degree structures. They first determined sufficient conditions for the existence of such trilayered graphs, then showed that if such subgraphs exist, a (well-placed) $\Theta$-graph could be embedded inside. Finally, they argued that either the preceding conditions hold, or the levels must expand by an average factor of $O(\frac{d}{\sqrt{k}\log k})$. Their result followed. 
 
In this paper, we follow the same proof structure. We improve on Bukh--Jiang's result by weakening the conditions required for minimum degree trilayered subgraphs to be present, and presenting a better method to embed well-placed $\Theta$-graphs in such subgraphs. These changes, presented in the following sections, lead to our $O(\sqrt{\log k})$ improvement on the best-known upper bound for $\ex(n, C_{2k})$. 

\section{Search for \texorpdfstring{$\Theta$}{Theta}-graphs}\label{sec:theta}
In this section, we present the central arguments of this paper. Our results are summarized in the following lemma, which states sufficient conditions for the existence of (well-placed) $\Theta$-graphs.

\begin{lemma}\label{lem:key}
Let $G$ be a trilayered graph with layers $V_1$, $V_2$, $V_3$, such that $d_{\text{min}}(G)\ge 2d + 5k^2$ and $d_{\text{max}}(G)\le \Delta d$. If the following conditions hold:
\begin{align}
d\cdot e(V_1, V_2)&\ge 40 k\log k\abs{V_3}, \label{c1} \\
e(V_1, V_2)&\ge 6k(\log k+1)^2(2\Delta k)^{2k-1}\abs{V_1}, \label{c2}\\
e(V_1, V_2) &\ge 20(\log k+1)\abs{V_2} \label{c3}
\end{align}
then either there is a $\Theta$-graph in $G[V_1, V_2]$, or 
there is a well-placed $\Theta$-graph in $G[V_1, V_2, V_3]$.
\end{lemma}

This lemma is an improvement over Lemma~6 in Bukh--Jiang. We removed two of the conditions and improved the last condition by a factor of $(t+1)$. 

To prove Lemma~\ref{lem:key}, the rest of this section is organized as follows: In Lemma~\ref{lem:base}, we show that given a trilayered graph formed by three consecutive levels in our BFS process, either we can find a trilayered subgraph with desired minimum degree structure, or we can find a trilayered subgraph with stronger constraints on its edges. This process can then be iterated --- In Lemma~\ref{lem:iterate}, we prove that under the conditions stated in Lemma~\ref{lem:key}, Lemma~\ref{lem:base} can be iterated to show the existence of a desired trilayered subgraphs. Finally, in Lemma~\ref{lem:embedding}, we show that a (well-placed) $\Theta$-graph can be embedded in such subgraphs, which completes the proof. 

Without further delay, we now quote the following result, which is Lemma~7 in Bukh--Jiang.
\begin{lemma}\label{lem:base}
Let $a,A,B,C,D$ be positive real numbers. Suppose $G$ is a trilayered graph with layers $V_1$, $V_2$, $V_3$ 
and the degree of every vertex in $V_2$ is at least $d+4k^2+C$. Assume also that
\begin{equation}\label{eq:base}
  a\cdot e(V_1,V_2)\ge (A+k+1)\abs{V_1}+B\abs{V_2}.
\end{equation}
Then one of the following holds:
\begin{enumerate}[label=\Roman{*}),ref=(\Roman{*})]
\item \label{alt:v1v2theta} There is a $\Theta$-graph in $G[V_1,V_2]$.

\item \label{alt:trilayer} There exist non-empty subsets $V_1'\subset V_1$, $V_2'\subset V_2$, $V_3'\subset V_3$ 
such that the induced trilayered subgraph $G[V_1',V_2',V_3']$ has minimum degree at least
$[A:B,C:D]$.

\item \label{alt:increment} There is a subset $\widetilde{V}_2\subset V_2$ such that $e(V_1,\widetilde{V}_2)\ge (1-a)e(V_1,V_2)$,
and $\abs{\widetilde{V}_2}\le D\abs{V_3}/d$.
\end{enumerate}
\end{lemma}

A proof of this lemma, as presented in Bukh--Jiang, is included in the Appendix for completeness.  Here the parameter $a$ can be interpreted as the edge loss ratio. More specifically, in case \ref{alt:v1v2theta} the proof of Lemma~\ref{lem:key} is complete, and similarly in case \ref{alt:trilayer} we are done by Lemma~\ref{lem:embedding}. In case \ref{alt:increment}, we found $\widetilde{V}_2$ that shrinks proportionally compared to $V_2$, while being adjacent to most edges between $V_1$ and $V_2$. We can then apply Lemma~\ref{lem:base} in $G[V_1, \widetilde{V}_2, V_3]$, thereby iterating this process. In the end, we will either obtain a subset of vertices of $V_2$ with an overly high average degree, or lands in case \ref{alt:v1v2theta} or \ref{alt:trilayer}. This procedure is done precisely in the following lemma, which is Lemma~8 in Bukh--Jiang. 

\begin{lemma}\label{lem:iterate}
Let $G$ be a trilayered graph with layers $V_1, V_2, V_3$ satisfying conditions \eqref{c1}, \eqref{c2} and \eqref{c3}. Let $C$ be a positive real number, such that the minimum degree from $V_2$ to $V_3$ is at least $d + 4k^2 + C$. Then one of the following holds:
\begin{enumerate}[label=\Roman{*}),ref=(\Roman{*})]
\item There is a $\Theta$-graph in $G[V_1, V_2]$. \label{outcome1}
\item There are non-empty subsets $V_1'\subset V_1$, $V_2'\subset V_2$, and $V_3'\subset V_3$ such that the induced trilayered subgraph $G[V_1', V_2', V_3']$ has minimum degree at least $[A:B, C:D]$, where $B\ge 5$, and
\begin{align}
A &\ge 2k(\Delta D)^{D-1}, \label{eq:const1} \\
(B-4)D &\ge 2k. \label{eq:const2}
\end{align}
\label{outcome2} 
\end{enumerate}
\end{lemma}

The following proof is an improved version of Bukh--Jiang's proof. 

\begin{proof}
Assume for the sake of contradiction that neither of the conclusions are true. We will first show that the conditions of Lemma~\ref{lem:base} hold for a tuple of well defined $A, B$ and $D$. Due to our assumptions, the only probable conclusion of Lemma~\ref{lem:base} would be \ref{alt:increment}, which gives us $\widetilde{V}_2\subseteq V_2$. We then iterate this procedure for $t = \log k$ steps on $\widetilde{V}_2$ and subsequent subsets of $V_2$. This process will generate a chain of sets $V_2^{(t)} \subseteq V_2^{(t-1)} \subseteq \cdots \subseteq V_2^{(1)}\subseteq V_2^{(0)} = V_2$. Finally, we will show a contradiction in $V_2^{(t)}$ to conclude the proof.

Let $a_i = \frac{1}{t-i+1}$, where $i$ ranges from $0$ to $t-1$. Let $V_2^{(0)} = V_2$. For $V_2^{(i)}$ that is well defined, set
\begin{align*}
d_i &= e(V_1, V_2^{(i)})/\abs{V_2^{(i)}}, \\
A_i &= a_ie(V_1, V_2^{(i)})/2\abs{V_1} - k - 1, \\
B_i &= a_id_i/4 + 5, \\
D_i &= \min(2k, 8k/a_id_i).
\end{align*}
Here $d_i$ is the average degree from $V_2^{(i)}$ to $V_1$. First note that $A_i, B_i, D_i$ satisfy constraints \eqref{eq:const1} and \eqref{eq:const2}. Indeed, \eqref{eq:const2} follows straightforwardly, and for \eqref{eq:const1}, we have by \eqref{c2}
\begin{align*}
A_i &= a_ie(V_1, V_2^{(i)})/2\abs{V_1} - k - 1 
\ge \frac{1}{2(t+1)^2}\frac{e(V_1, V_2)}{\abs{V_2}} - k - 1 \\
&\overset{\eqref{c2}}{\ge}3k(2\Delta k)^{2k-1}-k-1\ge 2k(\Delta D_i)^{D_i-1}.
\end{align*}
Therefore, if we apply Lemma~\ref{lem:base}, the only possible outcome is \ref{alt:increment}. The following claim is the key to our iteration process.
\begin{claim}
For $V_2^{(i)}$ that is well-defined, condition~\eqref{eq:base} of Lemma~\ref{lem:base} hold with respect to the above defined $a_i, A_i, B_i, C, D_i$. Moreover, let $V_2^{(i+1)}\subseteq V_2^{(i)}$ be the set derived from Lemma~\ref{lem:base}. We have the following invariants:
\begin{align}
e(V_1, V_2^{(i+1)}) &\ge (1-a_i)e(V_1, V_2^{(i)}), \label{eq:Invariant1}\\
d_{i+1} &\ge a_id_i\frac{t-i}{t+1}\frac{d\cdot e(V_1, V_2)}{8k\abs{V_3}}. \label{eq:Invariant2}
\end{align}
\end{claim}
\begin{proof}
This proof will proceed by induction. We first show that condition~\eqref{eq:base} holds for $i = 0$. 
\begin{align*}
(A_0+k+1)\abs{V_1} + B_i\abs{V_2} 
&= \frac{3}{4}a_0e(V_1, V_2) + 5\abs{V_2} \\
&\overset{\eqref{c3}}{\le} \frac{3}{4}a_0e(V_1, V_2) + \frac{1}{4(t+1)}e(V_1, V_2) = a_0e(V_1, V_2).
\end{align*}
Therefore, we can apply Lemma~\ref{lem:base} to obtain $\widetilde{V}_2\subseteq V_2$ as in outcome \ref{alt:increment}. Set $V_2^{(1)} = \widetilde{V}_2$. Invariant~\eqref{eq:Invariant1} then follows directly from the conclusions of Lemma~\ref{lem:base}. For \eqref{eq:Invariant2}, since $\abs{V_2^{(1)}} \le D_0\abs{V_3}/d$, we have
\begin{align*}
d_{1} &= \frac{e(V_1, V_2^{(1)})}{\abs{V_2^{(1)}}} \ge \frac{(1-a_0)e(V_1, V_2)}{D_0\abs{V_3}/d} \ge (1-a_0)a_0d_0\frac{d\cdot e(V_1, V_2)}{8k\abs{V_3}}.
\end{align*}
This completes the proof for the base case. For induction, note that iterative application of \eqref{eq:Invariant1} gives
\begin{align}
e(V_1, V_2^{(i)}) \ge e(V_1, V_2)\prod_{j=0}^{i-1}(1-a_j) = \frac{t-i+1}{t+1}e(V_1, V_2). \label{eq:Invariant3}
\end{align}
This inequality helps us show condition~\eqref{eq:base} again. Indeed, 
\begin{align*}
(A_i+k+1)\abs{V_1} + B_i\abs{V_2^{(i)}} 
&= \frac{3}{4}a_ie(V_1, V_2^{(i)}) + 5\abs{V_2^{(i)}} \le \frac{3}{4}a_ie(V_1, V_2^{(i)}) + 5\abs{V_2} \\
&\overset{\eqref{c3}}{\le} \frac{3}{4}a_ie(V_1, V_2^{(i)}) + \frac{1}{4(t+1)}e(V_1, V_2) \\
&\overset{\eqref{eq:Invariant3}}{\le} \frac{3}{4}a_ie(V_1, V_2^{(i)}) + \frac{t+1}{4(t+1)(t-i+1)}e(V_1, V_2^{(i)}) \\
&= \frac{3}{4}a_ie(V_1, V_2^{(i)}) + \frac{1}{4}a_ie(V_1, V_2^{(i)}) = a_ie(V_1, V_2^{(i)}).
\end{align*}
Therefore, by Lemma~\ref{lem:base} again, there is a subset $V_2^{(i+1)}\subset V_2^{(i)}$ satisfying \eqref{eq:Invariant1}, and
\[
\abs{V_2^{(i+1)}} \le D_i\abs{V_3}/d
\]
This implies
\begin{align*}
d_{i+1} &= \frac{e(V_1, V_2^{(i+1)})}{\abs{V_2^{(i+1)}}} \ge \frac{(1-a_i)e(V_1, V_2^{(i)})}{D_i\abs{V_3}/d} \ge (1-a_i)a_id_i\frac{d}{8k\abs{V_3}}e(V_1, V_2^{(i)}) \\
&\ge (1-a_i)a_id_i\frac{de(V_1, V_2)}{8k\abs{V_3}}\prod_{j=0}^{i-1}(1-a_j) \ge a_id_i\frac{t-i}{t+1}\frac{de(V_1, V_2)}{8k\abs{V_3}}.
\end{align*}
Therefore invariant~\eqref{eq:Invariant2} holds. This complete the proof of this claim.
\end{proof}
Through iterative application of this claim, we obtain our desired chain of subsets $V_2^{(t)} \subseteq \cdots \subseteq V_2^{(1)}\subseteq V_2^{(0)}$. For simplicity of notation, let $F = \frac{d\cdot e(V_1, V_2)}{8k\abs{V_3}}$. By \eqref{eq:Invariant2}, we have
\begin{align*}
d_i 
&\ge d_0 \cdot F^i\prod_{j=0}^{i-1}a_j\frac{t-j}{t+1} = d_0 \cdot F^i\prod_{j=0}^{i-1}\frac{t-j}{t-j+1}\frac{1}{t+1} \\
& = d_0 \cdot \Big(\frac{F}{t+1}\Big)^i\frac{t-i+1}{t+1} \overset{\eqref{c1}}{\ge} d_0 \cdot 5^i\Big(\frac{t}{t+1}\Big)^t\frac{t-i+1}{t+1} \\
&\ge d_0 \cdot 5^ie^{-1}\frac{t-i+1}{t+1}.
\end{align*}
Therefore we have
\begin{equation}
\frac{d_0}{d_i} \le \frac{e(t+1)}{5^i (t-i+1)}. \label{eq:deg}
\end{equation}
We now analyze the end results of our iteration, $V_2^{(t)}$ and $d_t$. Observe that $V_2^{(t)}$ preserves a good portion of the edges from $V_2$ to $V_1$ (invariant~\eqref{eq:Invariant1}), while having exponentially large average degree (equation~\eqref{eq:deg}). Therefore, we can draw different conclusions on the average degree of $G[V_1, V_2^{(t)}]$ dependent on $\abs{V_2^{(t)}}$. If we have $\abs{V_2^{(t)}}\le \abs{V_1}$, then 
\[
\frac{2e(V_1, V_2^{(t)})}{\abs{V_1} + \abs{V_2^{(t)}}} \ge \frac{e(V_1, V_2^{(t)})}{\abs{V_1}} \ge \frac{1}{t+1}\frac{e(V_1, V_2)}{\abs{V_1}}\ge 2k,
\]
which then implies outcome~\ref{outcome1} by Corollary~\ref{cor:avg_embed}. On the other hand, if $\abs{V_2^{(t)}}\ge \abs{V_1}$ and $d_t\ge 2k$, then 
\[
\frac{2e(V_1, V_2^{(t)})}{\abs{V_1} + \abs{V_2^{(t)}}} \ge \frac{e(V_1, V_2^{(t)})}{\abs{V_2^{(t)}}} = d_t\ge 2k,
\] 
which again leads to outcome~\ref{outcome1}. Therefore $d_t < 2k$ and
\[
d_0 \le d_te(t+1)/5^t < 2ke(t+1)/5^t < 20(t+1).
\]
This contradicts condition~\eqref{c1}. Therefore we conclude that the iteration must stops before $t$ steps, resulting in either outcome~\ref{outcome1} or outcome~\ref{outcome2}. \qedhere
\end{proof}

\begin{remark}\label{rem:one}
It could be shown in the above proof the contraction rate of $V_2^{(i)}$. Specifically,
\begin{align*}
\abs{V_2^{(i+1)}} &\le D_i\abs{V_3}/d \le \frac{1}{a_id_i}\frac{8k\abs{V_3}}{d} = \frac{d_0}{Fa_{i}d_i}\abs{V_2}\overset{\eqref{c1}}{\le} \frac{e(t+1)}{5^{i+1}t}\abs{V_2}.
\end{align*}
This bound confirms our intuition that $V_2^{(i)}$ shrinks exponentially. 
\end{remark} 

We now come to the last piece of the puzzle: proving the existence of a $\Theta$-graph. The following lemma, while following the same scheme as in Lemma~9 of Bukh--Jiang, presents a different method to embed an arbitrarily long path under the assumption that no (well-placed) $\Theta$-graphs exist. More details on such distinctions are discussed after the proof. 

\begin{lemma}\label{lem:embedding}
Let $G$ be a trilayered graph with layers $V_1, V_2, V_3$ and minimum degree at least $[A:B, d+k:D]$ where $B\ge 5$, and
\begin{align}
A &\ge 2k(\Delta D)^{D-1}, (B-4)D \ge 2k. \label{eq:para}
\end{align}
Assume that every vertex in $V_2$ has at most $\Delta d$ neighbors in $V_3$. Then there is a $\Theta$-graph in $G[V_2, V_3]$, or there is a well-placed $\Theta$-graph in $G$. 
\end{lemma}
\begin{proof}
Assume that neither of the conclusions are true. In this proof, we will utilize this assumption to embed an arbitrarily long path $P$ in $G$, contradicting the finiteness of the graph. $P$ will have the form $v_0\larrow v_1\larrow \dotsb \larrow v_l$, where $v_1, \cdots, v_l\in V_1$ and each pair $v_i, v_{i+1}$ is connected by a path of length $2D$ alternating between $V_2$ and $V_3$. 

To utilize the assumption of no well-placed $\Theta$-graph, we strengthen the statement by maintaining the following property while building the path:
\begin{definition}
A path $P$ is called good if every vertex in $V_2\cap P$ has at least one neighbor in $V_1\setminus P$.
\end{definition}
This property enables us to make arguments of the form ``either the path could be extended, or we can find a well-placed $\Theta$-graph'', as we will see later in the proof.

We start our construction with a random vertex $v_0$ from $V_1$. Inductively, assume that a good path $P=v_0\larrow v_1\larrow \dotsb \larrow v_{l-1}$ has been constructed, we wish to extend it to $v_0\larrow \cdots \larrow v_l$. We make the following observations.

\begin{claim}\label{cm:vi}
For all $i = 0, \cdots, l-1$, $v_i$ cannot have $k$ or more neighbors in $V_2\cap P$. 
\end{claim}
\begin{proof}
If $v_i$ has at least $k$ neighbors in $V_2\cap P$, then we can follow the path and build a $\Theta$-graph with a chord through $v_i$. This $\Theta$-graph is well-placed since $P$ is a good path.
\end{proof}

\begin{claim}\label{cm:small}
Given a good path $Q$, let $u\in V_2\cap Q$ be a vertex adjacent to the last vertex of $Q$ (note that this last vertex can belong to either $V_1$ or $V_3$). Then $u$ has less than $t = \lceil B/2\rceil$ neighbors in $V_1\cap Q$. 
\end{claim}
\begin{proof}
If $u$ has neighbors $v_{k_1},\cdots, v_{k_t},$ where $k_1 < k_2 < \cdots < k_t$, then the path $v_{k_2}\larrow u$ and the edge $uv_{k_2}$ form a cycle of length at least 
\[
2D(t-2) + 2\ge 2D(B/2 - 2) + 2 = D(B-4)+2\ge 2k.
\] 
This cycle, together with the chord $uv_{k_3}$, forms a $\Theta$-graph spanning over $V_1, V_2, V_3$. Moreover, this $\Theta$-graph is well-placed since $Q$ is a good path, and $u$ is adjacent to $v_{k_1}$ which is not part of the $\Theta$-graph. This contradicts our assumption. 
\end{proof}

Note that by Claim~\ref{cm:vi}, there are at least $A - k$ ways to extend $v_{l-1}$ to another vertex in $V_2\setminus P$, and Claim~\ref{cm:small} ensures that all of these extensions are good. Denote $U_0 = N(v_{l-1})\setminus P$, where $N(\cdot)$ is the usual notation for neighborhood. The following claim, which is the heart of our embedding scheme, states that a large portion of these good extensions in $U_0$ can be extended further inductively in a vertex-disjoint manner. 

\begin{claim}\label{cm:disjoint}
For $i = 0, 1, \cdots, D-1$, there exist sets $U_i\subset V_2$ such that for each $u\in U_i$, there exists a path $Q(u)$ from $U_0$ to $u$ of length $2i$ that alternates between $V_2$ and $V_3$. Moreover, $Q(u)$ is a good extension of $P$, and for every pair $u, v\in U_i$, $Q(u)$ and $Q(v)$ are vertex disjoint. Furthermore, 
\[
    \abs{U_i} \ge -3k + A \left(\frac{1}{8(2k+1)\Delta} \right)^i  \prod_{i=1}^{D-1}\frac{D-i}{i+1}.
\]  
\end{claim}
\begin{proof}
We prove this claim by induction, where the base case with $i = 0$ is true as stated. Assume the claim is true for $i$, we want to find $U_{i+1}$ by extending paths from $U_{i}$.  

For arbitrary $u\in U_i$, let $P_u$ denote the concatenation of paths $P$ and $Q(u)$. By similar argument as in Claim~\ref{cm:vi}, we see that $u$ cannot have more than $k$ neighbors in $P_u\cap V_3$. Therefore, $u$ has at least $d$ neighbors in $V_3$ that does not land on $P_u$. These neighbors are our candidates for extending $P_u$, and we filter these candidates with the following procedure. Define three sets $S_1, F_1$ and $T$, where $S_1, F_1\subset U_i$ and $T\subset V_3$. Intuitively, we want $S_1$ to be the set of vertices with successful extensions to $V_3$, $F_1$ to be $U_i\setminus S_1$, and $T$ to be the set of potential extensions from $S_1$ to $V_3$. Set them to be empty initually, consider the following procedure.
\begin{algorithm}\label{proc:one}
\caption{}
\begin{algorithmic}[1]
\State Pick a vertex $u$ randomly from $U_i\setminus (S_1\cup F_1)$.

\State Let $M_u = (N(u)\cap V_3) \setminus (T\cup P_u)$. If $\abs{M_u} \ge \frac{d}{2k+1}$, then randomly select $\frac{d}{2k+1}$ vertices in $M_u$ to put into $T$, and denote these vertices as $T_u$. Put $u$ into $S_1$. 

\State Otherwise, put $u\in F_1$ and move on to the next itration. Terminate this procedure if $S_1\cup F_1 = U_i$.
\end{algorithmic}
\end{algorithm}

We claim that when this procedure terminates, 
$
\abs{S_1}\ge \abs{U_i}/2.
$
Indeed, if $\abs{F_1} > \abs{U_i}/2$, then 
$
\abs{T} < \abs{S_1}\frac{d}{2k+1} < \frac{\abs{U_i}}{2}\frac{d}{2k+1}.
$
Moreover, every vertex $u$ in $F_1$ has at least $d$ neighbors in $V_3\setminus P_u$, which means at least $\frac{2kd}{2k+1}$ edges adjacent to $u$ land in $T$. Therefore, 
\[
e(F_1, T) \ge \abs{F_1}\frac{2kd}{2k+1} > \frac{2kd}{2k+1}\frac{\abs{U_i}}{2},
\]
which implies $e(F_1, T)/\abs{T} > 2k$. By Lemma~\ref{lem:min_embed}, there exists a $\Theta$-graph in $G[V_2, V_3]$, which is a contradiction. Thus $\abs{S_1} \ge \abs{U_i}/2$. 

We extend the previous notations to vertices in $T$. For $v\in T_u$ (as defined in Procedure~1), let $Q(v)$ be the path $Q(u)v$, and $P_v = P_uv$. Note that the paths $\{Q(v)\}_{v\in T}$ are not necessarily pairwise vertex disjoint, since $v$ could be on the path $Q(w)$ for some $w\in U_i, w\ne u$. This issue will be resolved later. For now, we make the following observation concerning extending vertices in $T$ back to $V_2$. 

\begin{claim}\label{cm:T2V2}
For an arbitrary vertex $v\in T$, it has at least $D - i$ neighbors in $V_2\setminus P_v$ or in the last $2k$ vertices of $P$. 
\end{claim}
\begin{proof}
For an edge $vw$ where $w\in P_v\cap V_2$, we call it long if the distance between $v, w$ is at least $2k$ through the path $P_v$ and short otherwise. If $v$ has a long edge $vw$, then $v$ cannot have any other neighbors in $P_v\cap V_2$, for otherwise there would be a well-placed $\Theta$-graph. Moreover, since $\abs{Q(v)\cap V_2} = i$, we see that $v$ has at most $i$ neighbors on $Q(v)$. The claim then follows. 
\end{proof}

Utilizing this claim, we will extend every vertex in $S_1$ greedily, while maintaining that all extensions land in different vertices in $V_2$. As in procedure~1, we define sets $S_2, F_2\subset S_1, D\subset V_2$, where $S_2$ denotes the set of vertices with successful 2-step extensions, and $F_2 = S_1\setminus S_2$. $D$ denotes the set of endpoints of successful extensions. We set them to be empty initially, and consider the following procedure.
\begin{algorithm}\label{proc:two}
\caption{}
\begin{algorithmic}[1]
\State Pick a vertex $u$ arbitrarily from $S_1\setminus (S_2\cup F_2)$.

\State If there exists $v\in N(u)\cap T$ and $w\in N(v)\setminus (P_v\cup D)$, then we can successfully extend $P_u$ to $P_uvw$. Put $u$ into $S_2$ and put $w$ into $D$.

\State If such vertices do not exist, put $u$ into $F_2$ and move on to the next itration. Terminate this procedure if $S_2\cup F_2 = S_1$.
\end{algorithmic}
\end{algorithm}

Let $\epsilon = \frac{D-i}{4(2k+1)\Delta}$. We claim that when this procedure terminates, $\abs{S_2} \ge \epsilon \abs{S_1} - 2k$. To see that, we know every vertex $u\in F_2$ cannot be extended, which means all of its possible extensions land in $D$ or the last $2k$ vertices of $P$. If $\abs{F_2} > (1-\epsilon)S_1 + 2k > (1-\epsilon)S_1$, by Procedure~1 and Claim~\ref{cm:T2V2}, the number of failed extension must be at least
\[
    \abs{F_2}\cdot \frac{d}{2k+1}\cdot (D-i) \ge \frac{(D-i)(1-\epsilon)d}{2(2k+1)}\abs{U_i}.
\]
Since $\abs{D} = \abs{S_2}$, all these failed extensions must land in a set of size less than $\epsilon \abs{S_1}$. The average degree on this set would then be at least
\[
     \frac{(D-i)(1-\epsilon)d}{2(2k+1)}\abs{U_i}\cdot \frac{1}{\epsilon \abs{U_i}} > 2(1-\epsilon)\Delta d > \Delta d,
\]
which is a contradiction to the assumption that no vertices in $V_2$ has more than $\Delta d$ neighbors in $V_3$. Therefore we have at least $\abs{S_2} \ge \epsilon\abs{U_i} - 2k$ successful extensions.

The next step is to filter these extensions such that they are pairwise vertex disjoint. 
What we have constructed so far is a set $Q$ of length $2i+2$ paths from $U_0$ to $D$ such that if we choose any two paths $p_1, p_2$ from $Q$, their first $2i$ vertices would be disjoint, and their last two vertices would also be disjoint. 
Therefore every path could only overlap with at most $2i+2$ other paths in $Q$, which implies there exists a set of pairwise disjoint paths $Q'$ such that $\abs{Q'} \ge \abs{Q}/(2i+2)$. 
Let $U_{i+1}\subset D$ be the set of endpoints of these paths, we have
\[
    \abs{U_{i+1}} = \abs{Q'} \ge \abs{Q}/(2i+2) = \abs{S_2}/(2i+2) \ge \frac{\epsilon}{2(i+1)}\abs{U_i} - 2k \ge \frac{D-i}{i+1}\frac{1}{8(2k+1)\Delta}\abs{U_i} - 2k,
\]
which satisfies the stated bound. All of these extensions are good by Claim~\ref{cm:small}. 
\end{proof}
Now from condition~\eqref{eq:para}, we see that $U_{D-1}$ is non-empty. Let $Q=v_0\larrow \dotsb \larrow v_{l-1} \larrow u$ be an arbitrary extension with $u\in U_{D-1}$. 
By Claim~\ref{cm:small}, $(N(u)\cap V_1)\setminus Q$ is non-empty. Let $v_l$ be chosen arbitrarily from this set, and let the new path be $Qv_l$. We prove one last claim to finish the proof.
\newpage
\begin{claim}
The path $Qv_l$ is good.
\end{claim}
\begin{proof}
We show that for any $w\in V_2\cap Q$, $w$ has at most $2t-2$ neighbors in $V_1\cap Qv_l$. By Claim~\ref{cm:small}, $w$ has fewer than $t$ neighbors in $Q\cap V_1$ that precede $w$ in~$Q$. We want to apply the same argument to the reverseal of $Qv_l$. 

Consider the sub-path $Q' = v_l\larrow w$ of $Q$. Since $Q$ is a good path, $w$ can't have $t$ or more neighbors in $V_1\cap Q$. Therefore, assume $w$ has neighbors $v_{k_1}, \cdots, v_{k_t}\in V_1\cap Q'$, where $v_{k_t} = v_l$ and $k_1 < k_2 < \cdots < k_t$. Then the path $v_{k_{t-1}} \larrow w$, together with the edges $wv_{k_{t-1}}$ forms a cycle of length at least $2k$, with chords through $v$. This $\Theta$-graph is well-placed since the path $Q$ is good, and $v_{k_{t-1}} \larrow w$ does not go through $v_{k_t} = v_l$, which means vertices of this $\Theta$-graph in $V_2$ can use $v_l$ to satisfy the well-placed condition. We conclude that $w$ must have less than $t$ neighbors in $Q'\cap V_1$. Since $2t-2<B$, the path $Qv_l$ is good.
\end{proof}
Therefore, we can construct an arbitrarily long path in $G$, which is a contradiction. We conclude that a (well-placed) $\Theta$-graph must exists.
\end{proof}
\begin{remark}\label{rem:two}
This result is stronger than Lemma~9 in Bukh--Jiang, in the sense that Bukh--Jiang showed how to embed one extension inductively, while we presented a method to embed multiple vertex-disjoint extensions simultaneously. We also note that $\abs{U_i}$ can be made arbitrarily large by increasing $A$, which only affects the magnitude of $n$. Therefore, our methods embed many ``parallel'' paths concurrently. 
\end{remark}

Our proof of Lemma~\ref{lem:key} is complete. We now proceed to prove Theorem~\ref{thm:real}. 

\section{Proof of Theorem~\ref{thm:real}}\label{sec:final}
In this section, we prove that under the conditions of Theorem~\ref{thm:real}, we have for all $i$
\begin{equation}
\abs{V_{i+1}} \ge (d^2/20k\log k)\abs{V_{i-1}}. \label{eq:expansion}
\end{equation}
We introduce the following auxiliary conditions, which will be proved by induction on $i$. 
\begin{align}
e(V_i, V_{i+1}) &\ge 2d\abs{V_{i}}, \label{eq:Bound1}\\
e(V_i, V_{i+1}) &\le 2k\abs{V_{i+1}}, \label{eq:Bound2}\\
\abs{V_{i+1}} &\ge k^{-1}d\abs{V_{i}},\label{eq:Bound3} 
\end{align}
These inequalities hold for $i = 0$. Assuming the inductive hypothesis, we know that the minimum degree in the graph is at least $2d+5k^2$. Therefore
\[
e(V_i, V_{i+1}) \ge (2d+5k^2)\abs{V_{i}} - e(V_{i-1}, V_i) \overset{\eqref{eq:Bound2}}{\ge} (2d+5k^2-2k)\abs{V_{i}} \ge 2d\abs{V_{i}}
\]
This inequality implies that $V_i$ has average degree at least $2d$ in $G[V_i, V_{i+1}]$. Moreover, if \eqref{eq:Bound2} is false, then $V_{i+1}$ has average degree at least $2k$ in $G[V_i, V_{i+1}]$. By Corollary~\ref{cor:avg_embed}, this leads to a contradiction. Therefore \eqref{eq:Bound2} is true, and \eqref{eq:Bound3} is a consequence of \eqref{eq:Bound1} and \eqref{eq:Bound2}. This completes the proof for the auxiliary claims.

We now move on to prove \eqref{eq:expansion}. Assume for the sake of contradiction that \eqref{eq:expansion} is false, we will show that the conditions of Lemma~\ref{lem:key} hold, which then leads to a contradiction.

Assume \eqref{c1} is false. We have
\begin{align*}
\frac{de(V_{i-1}, V_i)}{8k\abs{V_{i+1}}} &\le 5\log k, \\
2d^2\abs{V_{i-1}} \overset{\eqref{eq:Bound1}}{\le} de(V_{i-1}, V_i) &\le 40k\log k\abs{V_{i+1}}, \\
\abs{V_{i+1}} &\ge \frac{d^2}{20k\log k}\abs{V_{i-1}}.
\end{align*}
This contradicts with the assumption that \eqref{eq:expansion} is false. 

\eqref{c2} follows from the fact that $d \ge (20k)^{4k^2+2k}$. We have
\begin{align*}
6k(\log k + 1)^2 (2\Delta k)^{2k-1},
&\le 6k^3(2k^{3/2}(20k)^{2k}))^{2k-1} \\
&\le (20k)^{4k^2-2k}\cdot 6k^3\cdot (2k)^{3k} \le 2d \le e(V_i, V_{i+1})/\abs{V_i}.
\end{align*}
Finally, if \eqref{c3} is false, we have
\begin{align*}
2d\abs{V_{i-1}} &\overset{\eqref{eq:Bound1}}{\le} e(V_{i-1}, V_i) \le 20(\log k + 1)\abs{V_{i}} \overset{\eqref{eq:Bound3}}{\le} 40\log k\frac{k}{d}\abs{V_{i+1}}, \\
\abs{V_{i+1}} &\ge \frac{d^2}{20ek\log k}\abs{V_{i-1}}.
\end{align*}
This again implies \eqref{eq:expansion}. Therefore \eqref{eq:expansion} hold for all $i$. 

We now conclude the proof of Theorem~\ref{thm:real}. If $k$ is even, applying \eqref{eq:expansion} $k/2$ times results in
\[
\abs{V_{k}} \ge \frac{d^k}{(20k\log k)^{k/2}}.
\]
If $k$ is odd, applying \eqref{eq:expansion} $(k-1)/2$ times results in
\[
\abs{V_{k}} \ge \frac{d^{k-1}}{(20k\log k)^{(k-1)/2}} \abs{V_1} \ge \frac{d^{k}}{(20k\log k)^{k/2}}.
\]
Since $\abs{V_{k}} < n$, we must have $d < \sqrt{20k\log k}n^{1/k}$. \qed

\section{Conclusions}
We would like to point out that the constant factor of the upper bound proved by this paper is not fully optimized. In particular, one can further improve the bound by a factor of 2 if instead of using the Reduction Lemma, we employ a modified breadth-first search algorithm (see Bukh--Jiang Section~1) to bound the maximum degree in our graph.

\section{Acknowledgement}
The author would like to thank Boris Bukh for many constructive discussions on this problem, and for giving valuable comments on earlier versions of this paper. 

\section{Potential Improvements}\label{sec:imp}
This section is dedicated to devoted readers who intend to improve Theorem~\ref{thm:real} using our methods. 

The idea of using $\Theta$-graphs in the BFS approach originated from Pikhurko's work~\cite{Pikhurko}. The most critical component of this combination is the embedding scheme of a $\Theta$-graph in specific graph structures. For reference, Pikhurko utilized Lemma~\ref{lem:min_embed}, while Bukh--Jiang and the author utilized different versions of Lemma~\ref{lem:embedding}. In essence, all three proofs are driven by their respective embedding methods. Therefore, if one intends to improve the upper bound on $\ex(n, C_{2k})$ following this approach, one shall investigate potential structures and schemes to embed $\Theta$-graphs. 

We investigated the following structure in particular.
\begin{definition}
For $A, B, C, D\in \mathbb{R}$, we say that a trilayered graph $G$ on vertex sets $V_1, V_2, V_3$ has degree $[A: B, C: (2:D)]$ if there exists a partition of $V_2$ into $V_2^B$ and $V_2^C$, such that the minimum degree from $V_1$ to $V_2$, $V_2^B$ to $V_1$, $V_2$ to $V_3$, $V_3$ to $V_2^C$, and $V_3$ to $V_2^B$ are at least $A, B, C, 2, D$, respectively. 
\end{definition}

This definition is inspired by two observations. First of all, the proof of Lemma~\ref{lem:base} found that if an $[A:B, C:D]$ structure cannot be found, then $V_2$ can be partitioned into two sets $\widetilde{V}_2$ and $V_2\setminus \widetilde{V}_2$, such that the former has high density with $V_3$ and the latter has high density with $V_1$ (see Appendix). Let $V_C = \widetilde{V}_2$ and $V_B = V_2\setminus \widetilde{V}_2$, the existence of a $[A:B, C:(2,D)]$ structure is likely with respect to such graph partitions. Second, using the ideas in our proof of Lemma~\ref{lem:embedding} and Bukh--Jiang's proof of their Lemma~9, we can prove the following result.
\begin{lemma}
Under the same constraints on $A, B, D$ as in Lemma~\ref{lem:embedding}, if $G$ is a trylayered graph on $V_1, V_2 = V_B\cup V_C, V_3$ with minimum degree at least $[A:B, C:(2, D)]$, then there is a $\Theta$-graph in $G[V_2, V_3]$ or there is a well-placed $\Theta$-graph in $G$.
\end{lemma}
Therefore, if one is able to show, under weaker conditions in comparison to Lemma~\ref{lem:key}, that either an $[A:B, C:D]$ structure exists or an $[A:B, C:(2:D)]$ structure exists, then one could improve our bound. We were able to prove an analog of Lemma~\ref{lem:base} for the $[A:B, C:(2:D)]$ structure, but was unable to derive an analog of Lemma~\ref{lem:iterate}. 

\bibliographystyle{plain}
\bibliography{bib}

\section*{Appendix}
\paragraph{Bukh--Jiang's Proof of Reduction Lemma (slightly modified)}
Let $H$ be a subgraph of $G$ that maximizes the ratio $e(H)/v(H)^{1+\alpha/2}$. By the assumption on $e(G)$, this ratio is at least $cn^{\alpha/2}$.
Since $e(H)\le v(H)^2/2$, it then follows that $v(H)^{1-\alpha/2}\ge 2cn^{\alpha/2}$. 
Let $S$ be subset of $V(H)$ consisting of $\gamma v(H)$ vertices of largest degrees. We consider two cases.

Suppose at least $e(H)/4$ edges of $H$ are incident to vertices in $S$. Set $\eta=2\gamma/\alpha$. By averaging, we can find a set $T\subset V(H)\setminus S$ of $\eta v(H)$ elements that is incident
to at least fraction $\eta/(1-\gamma)$ of edges leaving $S$. Hence, $e(S\cup T)\ge (\frac{\eta}{1-\gamma})e(H)/4\ge \eta e(H)/4$. Let $H'$ be the subgraph of $H$ induced by $S\cup T$. 
Since
\begin{align*}
(\gamma+\eta)^{1+\alpha/2}
&=\gamma^{1+\alpha/2}(1+2/\alpha)^{1+\alpha/2} \le (3/\alpha)^{1+\alpha/2}\gamma^{1+\alpha/2} \\
&\le (3^{3/2}/\alpha^{1+\alpha/2})\gamma^{1+\alpha/2} \le (10/\alpha)\gamma^{1+\alpha/2} \le \gamma/2,
\end{align*}
we have
\[
    \frac{e(H')}{v(H')^{1+\alpha/2}} \ge \frac{\eta e(H)}{2\gamma v(H)^{1+\alpha/2}} = \frac{e(H)}{\alpha v(H)^{1+\alpha/2}} > \frac{e(H)}{v(H)^{1+\alpha/2}} ,
\] 
contradictory to the choice of $H$.

Therefore, we may assume that $S$ is incident to fewer than $e(H)/4$ edges of $H$. Thus the minimum degree of a vertex in $S$ is at most
$\frac{e(H)}{2\abs{S}} = \frac{e(H)}{2\gamma v(H)}$. Removing edges incident to $S$ from $H$ then leaves a graph $H'$ with maximum degree at most $ \frac{e(H)}{2\gamma v(H)}$ (since $S$ consists of vertices of highest degrees in $H$)
and at least $3e(H)/4$ edges. In particular, average degree of $H'$ is at least $3e(H)/(2v(H))$. 

Now we remove vertices of degree less than $e(H)/(2v(H))$ repeatedly to obtain $G'$. Since the number of edges removed is less than $e(H)/2$, $G'$ would have at least $\gamma v(H)$ vertices and $e(H)/4$ edges. Each vertex in this graph has degree between $e(H)/2v(H)$ and $e(H)/2\gamma v(H)$, and we have $e(G') \ge e(H)/4 \ge (c/4)n^{\alpha/2}v(H)^{1+\alpha/2}\ge (c/4)v(G)^{1+\alpha/2}$. Finally, since $e(H)/v(H)\ge cn^{\alpha/2}v(H)^{\alpha}\ge cv(G')^{\alpha}$, we are done.
\qed

\paragraph{Bukh--Jiang's Proof of Lemma~\ref{lem:base}} We suppose that alternative \ref{alt:v1v2theta} does not hold. Then, 
by Corollary~\ref{cor:avg_embed}, the average degree of 
every subgraph of $G[V_1,V_2]$ is at most $2k$.

Consider the process that aims to construct a subgraph satisfying \ref{alt:trilayer}.
The process starts with $V_1'=V_1$, $V_2'=V_2$ and $V_3'=V_3$,
and at each step removes one of the vertices that violate the minimum degree condition
on $G[V_1',V_2',V_3']$.
The process stops when either no vertices are left, or the minimum degree
of $G[V_1',V_2',V_3']$ is at least $[A:B,C:D]$. Since in the latter case we are done, we assume that
this process eventually removes every vertex of $G$. 

Let $R$ be the vertices of $V_2$ that were removed because at the time of removal
they had fewer than $C$ neighbors in $V_3'$. Put
\begin{align*}
  E'&\eqdef \{ uv\in E(G) : u\in V_2,\ v\in V_3,\ \text{and }v\text{ was removed before  }u\},\\
  S&\eqdef \{v \in V_2 : v\text{ has at least }4k^2\text{ neighbors in }V_1\}.
\end{align*}
Note that $\abs{E'}\le D\abs{V_3}$. We cannot have $\abs{S}\ge \abs{V_1}/k$, for
otherwise the average degree of the bipartite graph $G[V_1,S]$ would be at least $\frac{4k}{1+1/k}\ge 2k$. 
So $\abs{S}\le \abs{V_1}/k$.

The average degree condition on $G[V_1,S]$ implies that
\[
  e(V_1,S)\le k(\abs{V_1}+\abs{S})\le (k+1)\abs{V_1}.
\]

Let $u$ be any vertex in $R\setminus S$. Since it is connected to
at least $(d+4k^2+C)-4k^2 = d+C$ vertices of $V_3$, it must be adjacent to at least $d$ 
edges of $E'$. Thus,
\[
  \abs{R\setminus S}\le \abs{E'}/d\le D\abs{V_3}/d.
\]

Assume that the conclusion \ref{alt:increment} does not hold with $\widetilde{V}_2=R\setminus S$. Then
$e(V_1,R\setminus S)<(1-a)e(V_1,V_2)$. Since the total
number of edges between $V_1$ and $V_2$ that were removed
due to the minimal degree conditions on $V_1$ and $V_2$ 
is at most
$A\abs{V_1}$ and $B\abs{V_2}$ respectively, we conclude that
\begin{eqnarray*}
  e(V_1,V_2)&\le& e(V_1,S)+e(V_1,R\setminus S)+ A\abs{V_1}+B\abs{V_2}\\
            &<& (k+1)\abs{V_1}+(1-a)e(V_1,V_2)+A\abs{V_1}+B\abs{V_2},
\end{eqnarray*}
implying that
$$a\cdot e(V_1,V_2) < (A+k+1)\abs{V_1}+B\abs{V_2}.$$
The contradiction with \eqref{eq:base} completes the proof. \qed

\end{document}